\documentclass[reqno, 12pt]{amsart}
\title[Recurrence or transience of random walks on random graphs]{Recurrence or transience of random walks on random graphs generated by point processes in $\RR^d$}
\author[A. Rousselle]{Arnaud Rousselle$^\ddagger$}
\thanks{$\ddagger$ Normandie Universit\'e, Universit\'e de Rouen, Laboratoire de Math\'ematiques Rapha\"el Salem, CNRS, UMR 6085, Avenue de l'universit\'e, BP 12, 76801 Saint-Etienne du Rouvray Cedex, France.
\newline E-mail address: arnaud.rousselle1@univ-rouen.fr}
\usepackage{latexsym,amssymb,amsmath,a4wide,hyperref,alltt,cite, graphicx,enumerate}

\def\d{\mathrm{d}}

\def\x{\mathbf{x}}
\def\y{\mathbf{y}}
\def\z{\mathbf{z}}

\def\N{\mathcal{N}}

\def\PP{\mathbb{P}}
\def\RR{\mathbb{R}}
\def\EE{\mathbb{E}}
\def\NN{\mathbb{N}}
\def\ZZ{\mathbb{Z}}
\def\DT{\operatorname{DT}}
\def\Vor{\operatorname{Vor}}

\def\II{\mathrm{ 1\hskip -.29em I}}

\newtheorem*{ass*}{Assumptions}
\newtheorem{defi}{Definition}
\newtheorem{lemm}[defi]{Lemma}
\newtheorem{crit}[defi]{Criterion}
\newtheorem{theo}[defi]{Theorem}
\newtheorem{rem}[defi]{Remark}

\newenvironment{dem}{\vskip 2mm\noindent {\it Proof}:}
                    {\hfill $\square$ \vskip 2mm \noindent}

\begin{document}
\begin{abstract}
We consider random walks associated with conductances on Delaunay triangulations, Gabriel graphs and skeletons of Voronoi tilings generated by point processes in $\RR^d$. Under suitable assumptions on point processes and conductances, we show that, for almost any realization of the point process, these random walks are recurrent if $d=2$ and transient if $d\geq 3$. These results hold for a large variety of point processes including Poisson point processes,  Mat\'ern cluster and Mat\'ern hardcore processes which have clustering or repulsive properties. In order to prove them, we state general criteria for recurrence or transience which apply to random graphs embedded in $\RR^d$.  

\noindent \emph{To appear in Stochastic Processes and their Applications.}
\end{abstract}
\date{\today}
\maketitle

{\bf Key words:} Random walk in random environment; recurrence; transience; Voronoi tessellation; Delaunay triangulation; Gabriel graph; point process; electrical network. 

{\bf AMS 2010 Subject Classification :} Primary: 60K37, 60D05; secondary: 60G55; 05C81.

\section{Introduction and main results}
We deal with the question of recurrence or transience for random walks on Delaunay triangulations, Gabriel graphs and skeletons of Voronoi tilings which are generated by point processes in $\RR^d$. We obtain recurrence or transience results for random walks on these graphs for a.a. realization of the underlying point process. Random walks on random geometric networks are natural for describing flows, molecular diffusions, heat conduction or other problems from statistical mechanics in random and irregular media. In the last decades, most authors considered models in which the underlying graphs were the lattice $\ZZ^d$ or a subgraph of $\ZZ^d$. Particular examples include random walks on percolation clusters (see \cite{Grimmett93} for the question of recurrence or transience and \cite{BergerBiskup,MathieuQuenched,MP} for invariance principles) and the so-called random conductance model (see \cite{BiskupLN} and references therein). In recent years, techniques were developed to analyze random walks on complete graphs generated by point processes with jump probability which is a decreasing function of the distance between points. The main results on this model include an annealed invariance principle \cite{FSS,FM}, isoperimetric inequalities \cite{CF}, recurrence or transience results \cite{CFG}, a quenched invariance principle and heat kernel estimates \cite{CFP}. If the underlying graph is not the complete graph but, for example, a Delaunay triangulation or a Gabriel graph, additional difficulties appear due to the lack of information on the structure of the graphs (distribution and correlations of the degrees or of the edge lengths, volume growth, \dots). Very recently, Ferrari, Grisi and Groisman \cite{FGG} constructed harmonic deformations of Delaunay triangulations generated by point processes which could be used to prove a quenched invariance principle via the corrector method, at least in the 2-dimensional case. Unfortunately, the sublinearity of the corrector does not directly follow from their proofs in higher dimensions. Precise heat kernel estimates, similar to those derived by Barlow \cite{Barlow} in the percolation setting, are still to be obtained.

Under suitable assumptions on the point process and the transition probabilities, we obtain in this paper almost sure recurrence or transience results, namely Theorem \ref{thprinc} and Theorem \ref{proprinc}, for random walks on the Delaunay triangulation, the Gabriel graph or the skeleton of the Voronoi tiling of this point process. We hope that the different methods developed in the papers \cite{FSS,FM,CF,CFP} can be adapted to our setting to obtain annealed and quenched invariance principles. This will be the subject of future work. 

\subsection{Conditions on the point process}\label{condproc} 
In what follows the point process $\N$ is supposed to be simple, stationary and almost surely in general position (see \cite{Zessin}): a.s. there are no $d+1$ points (resp. $d+2$ points) in any $(d-1)$-dimensional affine subspace (resp. in a sphere).

In this paper, $c_1, c_2,\dots$ denote positive and finite constants. We will need the following assumptions on the void probabilities {\bf (V)} and on the deviation probabilities {\bf (D$_2$)} and {\bf (D$_{3^+}$)}:
\begin{ass*}$\,$
\begin{itemize}
\item[{\bf (V)}]
There exists a constant $c_1$ such that for $L$ large enough:
\begin{equation*}
\PP\big[\#\big([0,L]^d\cap\N\big)=0\big]\leq e^{-c_1L^d}.
\end{equation*}
\item[{\bf (D$_2$)}] If $d=2$, there are constants $c_2,c_3$ such that for $L,l$ large enough:
\begin{equation*}
\PP\big[\#\big(\big([0,L]\times[0,l]\big)\cap\N\big)\geq c_2Ll\big]\leq e^{-c_3Ll}.
\end{equation*}
\item[{\bf (D$_{3^+}$)}] \label{condproc3} If $d\geq 3$, there exists $c_4$ such that for $L$ large enough and all $m>0$:
\begin{equation*}
\PP\big[\#\big([0,L]^d\cap\N\big)\geq m\big]\leq e^{c_4L^d-m}.
\end{equation*}
\end{itemize}
\end{ass*}

For transience results, the following additional assumptions are needed:
\begin{itemize}
\item[{\bf (FR$_k$)}] \emph{$\N$ has a \emph{finite range of dependence $k$}, i.e., for any disjoint Borel sets $A, B\subset \RR^d$ with $\operatorname{d}(A,B):=\inf\{\Vert x-y\Vert\,:\,x\in A,y\in B\}\geq k$, $\N\cap A$ and $\N\cap B$ are independent.} 
\item[{\bf (ND)}] \emph{Almost surely, $\N$ does not have any \emph{descending chain}, where a descending chain is a sequence $(u_i)_{i\in\NN}\subset\N$ such that: $$\forall i\in\NN,\,\Vert u_{i+2}-u_{i+1}\Vert <\Vert u_{i+1}-u_i\Vert.$$}
\end{itemize}
As discussed in Section \ref{exproc}, these assumptions are in particular satisfied if $\N$ is:
\begin{itemize}
\item a homogeneous Poisson point process (PPP),
\item a Mat\'ern cluster process (MCP),
\item a Mat\'ern hardcore process I or II (MHP I/II).
\end{itemize}
Moreover, assumptions {\bf (V)} and {\bf (D$_2$)} hold if $\N$ is a stationary determinantal point process (DPP). 

A brief overview on each of these point processes is given in Section \ref{exproc}. Note that these processes have different interaction properties: for PPPs there is no interaction between points, MCPs exhibit clustering effects whereas points in MHPs and DPPs repel each other. 
\subsection{The graph structures}
Write $\mathrm{Vor}_\N(\x):=\{x\in\RR^d\, :\,\Vert x-\x\Vert\leq\Vert x-\y\Vert,\forall\y\in\N \}$ for the Voronoi cell of $\x\in\N$; $\x$ is called the nucleus or the seed of the cell. The Voronoi diagram of $\N$ is the collection of the Voronoi cells. It tessellates $\RR^d$ into convex polyhedra. See \cite{Moller,NewPerspectives} for an overview of these tessellations.

The graphs considered in the sequel are:  
\begin{description}
\item[{\bf $\operatorname{VS}(\N)$}]
the skeleton of the Voronoi tiling of $\N$. Its vertex (resp. edge) set consists of the collection of the vertices (resp. edges) on the boundaries of the Voronoi cells. Note that this is the only graph with bounded degree considered in the sequel. Actually, if $\N$ is in general position in $\RR^d$,  any vertex of $\operatorname{VS}(\N)$ has degree $d+1$. 
\item[{\bf $\operatorname{DT}(\N)$}] the Delaunay triangulation of $\N$. It is the dual graph of its Voronoi tiling of $\mathcal{N}$. It has $\N$ as vertex set and there is an edge between $\x$ and $\y$ in $\operatorname{DT}(\N)$ if $\mathrm{Vor}_\N(\x)$ and $\mathrm{Vor}_\N(\y)$ share a $(d-1)$-dimensional face. Another useful characterization of $\operatorname{DT}(\N)$ is the following: a simplex $\Delta$ is a cell of $\operatorname{DT}(\N)$ \emph{iff} its circumscribed sphere has no point of $\N$ in its interior. Note that this triangulation is well defined since $\N$ is assumed to be in general position. 

These two graphs are widely used in many fields such as astrophysics \cite{Ramella}, cellular biology \cite{Poupon}, ecology \cite{Roque} and telecommunications \cite{BB3}. 
\item[{\bf $\operatorname{Gab}(\N)$}] the Gabriel graph of $\N$. Its vertex set is $\N$ and there is an edge between $u,v\in \N$ if the ball of diameter $[u,v]$ contains no point of $\N$ in its interior. Note that $\operatorname{Gab}(\N)$ is a subgraph of $\operatorname{DT}(\N)$ (see \cite{Matula}). It contains the Euclidean minimum spanning forest of $\N$ (in which there is an edge between two points of $\N$ $\x$ and $\y$ \emph{iff} there do not exist an integer $m$ and vertices $u_0=\x, \dots, u_m=\y\in\N$ such that $\Vert u_i-u_{i+1}\Vert <\Vert \x-\y\Vert$ for all $i\in\{0, \dots, m-1\}$, see \cite{AS}) and the relative neighborhood graph (in which there is an edge between two points $x$ and $y$ of $\N$ whenever there does not exist a third point that is closer to both $x$ and $y$ than they are to each other). It has for example applications in geography, routing strategies, biology and tumor growth analysis (see\cite{BBD,Gab1,Gab2}).
\end{description}

\subsection{Conductance function}
Given a realization $\N$ of a point process and an unoriented graph $G(\N)=(V_{G(\N)},E_{G(\N)})$ obtained from $\N$ by one of the above constructions,  a \emph{conductance} is a positive function on $E_{G(\N)}$. For any vertex $u$ set $w(u):=\sum_{v\sim u}C(u,v)$ and $R=1/C$ the associated \emph{resistance}. Then $({G(\N)},C)$ is an infinite \emph{electrical network} and the random walk on ${G(\N)}$ associated with $C$ is the (time homogeneous) Markov chain $(X_n)_n$ with transition probabilities given by:
\[\PP\big[X_{n+1}=v\big| X_n=u\big]=\frac{C(u,v)}{w(u)}.\] 
Note that the random walk is well defined for locally finite graphs and that the graphs of the three classes considered
in the previous section are a.s. locally finite.
We refer to \cite{DoyleSnell,Lyons} for introductions to electrical networks and random walks associated with conductances. In our context, the cases where $C$ is either constant on the edge set (simple random walk on ${G(\N)}$) or given by a decreasing positive function of edge length are of particular interest. Our models fit into the broader context of random walks on random graphs with conductances, but due to the geometry of the problem we need to adapt the existing techniques. 

\subsection{Main results}   
The obtained random graph is now equipped with the conductance function $C$ possibly depending on the graph structure. Our main result is:  
\begin{theo}\label{thprinc} Let $\N$ be a homogeneous Poisson point process, a Mat\'ern cluster process or a Mat\'ern hardcore process of type I or II.
\begin{enumerate} 
\item Let $d=2$. Assume that, for almost any realization of $\N$, $C$ is bounded from above by a finite constant whose value possibly depends on $\N$. Then for almost any realization of $\N$ the random walks on $\operatorname{DT}(\N)$, $\operatorname{Gab}(\N)$ and $\operatorname{VS}(\N)$ associated with $C$ are recurrent.
\item Let $d\geq 3$. If $C$ is uniformly bounded from below or a decreasing positive function of the edge length, then for almost any realization of $\N$ the random walks on $\operatorname{DT}(\N)$, $\operatorname{Gab}(\N)$ and $\operatorname{VS}(\N)$ associated with $C$ are transient.
\end{enumerate}

Moreover, \emph{(1)} holds  if $\N$ is a stationary determinantal point process in $\RR^2$.
\end{theo}
Note that, in $(1)$, the bounded conductances case follows immediately from the unweighted case (i.e. $C\equiv1$) by Rayleigh monotonicity principle and we can restrict our attention to this last case in the proofs.

To the best of our knowledge, recurrence or transience of random walks on this kind of graphs has been sparsely considered
in the literature. Only in the unpublished manuscript \cite{AddarioberrySarkar05} Addario-Berry and Sarkar announced similar results in the setting of simple random walks on the Delaunay triangulation generated by a PPP and they noticed that their method can be applied to more general point processes. Their proofs relied on a deviation result for the so-called \emph{stabbing number} of $\operatorname{DT}(\N)$ contained in a second unpublished manuscript (\emph{The slicing number of a Delaunay triangulation} by Addario-Berry, Broutin, and Devroye). This last work is unfortunately unavailable. Note that the deviation result for the \emph{stabbing number} has been proved since then in \cite{Pimentel08}. We develop a new method, which avoids the use of such a strong estimate and is thus more tractable. This allows us to obtain recurrence and transience results for a large class of point processes and geometric graphs. 

Besides, several works show that random walks on distributional limits of finite rooted planar random graph are almost surely recurrent (see Benjamini and Schramm \cite{BenjaminiSchramm} and Gurel-Gurevich and Nachmias \cite{GurelGurevichNachmias}). Let us briefly explain how the results of \cite{GurelGurevichNachmias} could be used to obtain the recurrence of simple random walks on Delaunay triangulations generated by Palm measures associated with point processes in the plane.  Given $\N$ and $n$, consider the sub-graph $G_n(\N)$ of $\DT(\N)$ with vertex set given by $V_n(\N):=\{\x \in\N: \Vor_\N(\x)\cap [-n,n]^2\neq \emptyset\}$ and with edge set induced by the edges of $\DT(\N)$. If the point process is stationary, one could obtain Delaunay triangulations generated by the Palm version of point processes rooted at 0 as distributional limits of $(G_n,\rho_n)$ where, for each $n$, $\rho_n$ is chosen uniformly at random in $V_n(\N)$. Thanks to the results of Zuyev \cite{Zuyev}, we know that the degree of the origin in the Delaunay triangulation generated by the Palm measure of a PPP has an exponential tail. The proof can be adapted in the non-Poissonian case under assumptions similar to {\bf (V)} and {\bf (D$_2$)} for the Palm version of the point process and an assumption of finite range of dependence. Using a similar construction, one could expect to apply the results of \cite{BenjaminiSchramm} to the skeleton of the Voronoi tiling $\operatorname{VS(\N)}$ in the plane in which each vertex has degree 3. In this case, one must consider the Palm measure of the point process of the vertices of $\operatorname{VS(\N)}$ on which there is no information. The main problem of this approach is that it provides a recurrence result for random walks on the Palm version of the point process of $\operatorname{VS}$-vertices. It has no clear connection with the initial point process of the nuclei of the Voronoi cells. We do think that the recurrence criterion stated below is well adapted to the particular geometric graphs that we consider in this paper.

In the sequel, we will actually prove Theorem \ref{proprinc} which implies Theorem \ref{thprinc} and deals with general point processes as described in Subsection \ref{condproc}: 
\begin{theo}\label{proprinc}
Let $\N$ be a stationary simple point process in $\RR^d$ almost surely in general position.  
\begin{enumerate}
\item \label{proprincD2} Let $d=2$. Assume that $\N$ satisfies {\bf (V)} and {\bf (D$_2$)}. If $C$ is uniformly bounded from above, then for almost any realization of $\N$ the random walks on $\operatorname{DT}(\N)$, $\operatorname{Gab}(\N)$ and $\operatorname{VS}(\N)$ associated with $C$ are recurrent.
\item Let $d\geq 3$. Assume that $\N$ satisfies {\bf (V)}, {\bf (D$_{3^+}$)} and {\bf (FR$_k$)}. If $C$ is uniformly bounded from below or a decreasing positive function of the edge length, then for almost any realization of $\N$ the random walks on $\operatorname{DT}(\N)$ and $\operatorname{VS}(\N)$ associated with $C$ are transient.

If in addition $\N$ satisfies {\bf (ND)}, the same conclusion holds on $\operatorname{Gab}(\N)$.
\end{enumerate}
\end{theo}
\begin{rem}\label{remNonStat} The stationarity assumption is not required. Indeed, one can derive similar results when the underlying point process is not stationary and {\bf (V)}, {\bf (D$_2$)} and {\bf (D$_{3^+}$)} are replaced by:
\begin{itemize}
\item[{\bf (V')}]
There exists a constant $c'_1$ such that for $L$ large enough:
\begin{equation*}
\PP\big[\#\big(\big( a+[0,L]^d\big)\cap\N\big)=0\big]\leq e^{-c'_1L^d},\qquad \forall a\in\RR^d.
\end{equation*}
\item[{\bf (D'$_2$)}] If $d=2$, there are constants $c'_2,c'_3$ such that for $L,l$ large enough:
\begin{equation*}
\PP\big[\#\big(\big(a+\big([0,L]\times[0,l]\big)\big)\cap\N\big)\geq c'_2Ll\big]\leq e^{-c'_3Ll},\qquad \forall a\in\RR^d.
\end{equation*}
\item[{\bf (D'$_{3^+}$)}] If $d\geq 3$, there exists $c'_4$ such that for $L$ large enough and all $m>0$:
\begin{equation*}
\PP\big[\#\big(\big(a+[0,L]^d\big)\cap\N\big)\geq m\big]\leq e^{c'_4L^d-m},\qquad \forall a\in\RR^d.
\end{equation*}
\end{itemize}

Let us also point out that conditions {\bf (V)}, {\bf (D$_{3^+}$)} and {\bf (FR$_k$)} can be replaced by the domination assumption (\ref{condtrans3}) of Criterion \ref{crittrans} for the processes of \emph{good boxes} defined in Subsections \ref{GBVor} and \ref{GBGab}. In particular, this allows to relax the finite range of dependence assumption and to obtain results for point processes with good mixing properties (see also Remark \ref{remtrans1}).
\end{rem}

\subsection{Outline of the paper}
The use of the theory of electrical networks allows us to derive recurrence and transience criteria, namely Criteria \ref{critrec} and \ref{crittrans}, which are well suited to the study of random walks on random geometric graphs embedded in the Euclidean space. These general criteria are proved in a concise way in Sections \ref{secrec} and \ref{sectrans}.  Sections \ref{dim2} and \ref{dim3} are devoted to the proof of Theorem \ref{proprinc}. They constitute the heart of the paper and rely on arguments from stochastic geometry. In Section \ref{exproc}, assumptions of Theorem \ref{proprinc} are proved to hold for the point processes considered in Theorem \ref{thprinc}.
\section{A recurrence criterion}\label{secrec}

In this section, we give a recurrence criterion for random walks on a graph $G=(V_G,E_G)$ embedded in $\RR^d$ and equipped with a conductance $C$. In the sequel $G$ is assumed to be connected, infinite and locally finite. The criterion is established on deterministic graphs, but it gives a way to obtain almost sure results in the setting of random graphs (see Section \ref{dim2}). Note that it is a slight generalisation of the proof of \cite[Theorem 4]{AddarioberrySarkar05} and is close in spirit to the Nash-Williams criterion (see \cite[\S 2.5]{Lyons}). 

Let us define:
\begin{align*}
A_i&:=\{x\in \RR^d:\,i-1\leq \Vert x\Vert_\infty <i\},\\
B_i&:=[-i,i]^d,\\
Ed_G(i)&:=\{e=(u,v)\in E_G: u\in B_i,v \not\in B_i\}.
\end{align*}

\begin{crit}\label{critrec}
Assume that the conductance $C$ is bounded from above and there exist functions $L$ and $N$ defined on the set $\NN^*$ of positive integers and with values in the set $\RR^*_+$ of positive real numbers such that:
\begin{enumerate}
\item \label{cond1}$\underset{i}{\sum}\frac{1}{L(i)N(i)}=+\infty$, 
\item \label{cond2}the length of any edge in $Ed_G(i)$ is less than $L(i)$,
\item \label{cond3}the number of edges in $Ed_G(i)$ is less than $N(i)$.
\end{enumerate}

Then the random walk on $G$ with conductance $C$ is recurrent.
\end{crit}
\begin{rem}
\begin{enumerate}
\item One can check that $w(u)=\sum_{v\sim u}C(u,v),\,u\in V_G$, is a reversible measure of infinite mass, hence the random walk is \emph{null recurrent}.
\item The estimates appearing in the proof of this criterion lead to lower bounds for the effective resistance between the closest point $x_0$ to the origin and $\N\cap B_n^c$. This quantity is of interest because it is related to the expected number of visits of the random walk at $x_0$ before leaving $B_n$. For each example given in this paper, we could obtain lower bounds to be of order $\log(\log (n))$; this is not expected to be the correct order. It seems credible that the correct order is, as for simple walks on $\ZZ^2$, $\log (n)$.
\end{enumerate}
\end{rem}
\begin{dem}
The basic idea is inherited from electrical network theory. One can define the \emph{effective resistance to infinity} of the network which is known to be  infinite if and only if the associated random walk is recurrent (see \cite[\S 2.2]{Lyons}). We will reduce the network so that the resistance of any edge does not increase. Thanks to the Rayleigh monotonicity principle, the effective resistance to infinity does not increase either. It then suffices to show that the reduced network has infinite resistance to infinity. 

First each edge $e=(u,v)$ with $u\in A_{i_1},\, v\in A_{i_2}$, $i_1<i_2$ is cut into $j=i_2-i_1$ resistors connected in series with endpoints in consecutive annuli $A_{i_1},\dots ,A_{i_2}$ each one having resistance $(jC(e))^{-1}$. Secondly, points of a annulus $A_i$ are merged together into a single point $a_i$. We hence obtain a new network with vertices $(a_i)_{i\geq i_0}$ where $i_0$ is the lowest index such that $A_{i_0}$ contains a vertex of the original graph $G$. The resistance to infinity of the new network is lower than the original one and is equal to $\sum_{i=i_0}^{+\infty}r_i$, where $r_i$ stands for the effective resistance between $a_i$ and $a_{i+1}$.

It remains to show that:
\[\sum_{i=i_0}^{+\infty}r_i=+\infty.\]

Let us subdivide $Ed_G(i)$ into the following subsets:
\[Ed_G(i,j):=\big\{e=(u,v)\in E_G:\,u\in A_{i_1},\, v\in A_{i_2},\, i_1\leq i< i_2,\, i_2-i_1=j\big\}.\]

Note that an edge $e\in Ed_G(i,j)$ provides a conductance of $jC(e)$ between $a_i$ and $a_{i+1}$ in the new network. Thanks to the usual reduction rules and conditions (\ref{cond2}) and (\ref{cond3}), we get:
\begin{align*}
\frac{1}{r_i}&=\sum_{j=1}^{+\infty}\sum_{e\in Ed_G(i,j)}jC(e)
\leq (\sup C)\sum_{j=1}^{+\infty}j\#Ed_G(i,j)\\
&\leq (\sup C)\sum_{j=1}^{L(i)}j\#Ed_G(i,j)
\leq (\sup C) L(i)\sum_{j=1}^{L(i)}\#Ed_G(i,j)\\
&\leq (\sup C)L(i)\sum_{j=1}^{+\infty}\#Ed_G(i,j)
\leq (\sup C) L(i)N(i).
\end{align*}
This concludes the proof since:
\[\sum_{i=i_0}^{+\infty}r_i\geq \sum_{i=i_0}^{+\infty}\frac{1}{(\sup C)L(i)N(i)}\]
and the r.h.s. is infinite  by condition (\ref{cond1}).
\end{dem}
\section{A transience criterion}\label{sectrans}
In this section, the graph $G=(V_G,E_G)$ is obtained from a point process in $\RR^d$ and is equipped with a conductance $C$ (and therefore a resistance $R$) on $E_G$. The key idea is to combine \emph{discretization techniques} with a \emph{rough embedding method}. 

Let us define \emph{rough embeddings} for unoriented networks in a similar way to \cite[\S 2.6]{Lyons}:  

\begin{defi} \label{DefRoughEmb}
Let $H$ and $H'$ be two networks with resistances $r$ and $r'$.

We say that a map $\phi : V_H \longrightarrow V_{H'}$ is a \emph{rough embedding from $(H,r)$ to $(H',r')$} if there exist $\alpha, \beta <+\infty$ and a map $\Phi$ from (unoriented) edges of $H$ to (unoriented) paths in $H'$ such that:
\begin{enumerate}
\item for every edge $(u, v)\in E_H$, $\Phi (u, v)$ is a non-empty simple ({\it i.e.} with no repeating vertices) path of edges of $H'$ between $\phi (u)$ and $\phi (v)$ with
\[\sum_{e'\in\Phi(u, v)}r'(e')\leq \alpha r(u,v);\]
\item any edge $e'\in E_{H'}$  is in the image under $\Phi$ of at most $\beta$ edges of $H$.
\end{enumerate}
\end{defi}

A version of \cite[Theorem 2.17]{Lyons} without orientation is needed to establish the criterion. Lyons and Peres attribute this result to Kanai \cite{Kanai86}. 
\begin{theo}[see \cite{Lyons}]\label{roughembbeding}
If there is a rough embedding from $(H,r)$ to $(H',r')$ and $(H,r)$ is transient, then $(H',r')$ is transient. 
\end{theo}
To deduce Theorem \ref{roughembbeding} from \cite[Theorem 2.17]{Lyons}, one can do the following. We construct oriented graphs $\overrightarrow{H}$ and $\overrightarrow{H}'$ from $H$ and $H'$ such that if there is an unoriented edge between $x$ and $y$ in $H$ (resp. in $H'$) with resistance $r(\{x,y\})$, there are oriented edges from $x$ to $y$ and from $y$ to $x$ in $\overrightarrow{H}$ (resp. in $\overrightarrow{H}'$) with resistance $r(\{x,y\})$. It is then easy to see that if there is a rough embedding from $(H,r)$ to $(H',r')$ in the sense of Definition \ref{DefRoughEmb}, there is also a rough embedding from $(\overrightarrow{H},r)$ to $(\overrightarrow{H'},r')$ according to the definition given in \cite{Lyons}. 

Let us divide $\RR^d$ into boxes of side $M\geq 1$: \[B_\z=B^M_\z:=M\z+\Big[-\frac{M}{2},\frac{M}{2}\Big)^d,\,\z\in \ZZ^d.\] 

We will prove the following criterion:
\begin{crit}\label{crittrans}
If $d\geq 3$ and if one can find a subset of boxes, called \emph{good boxes}, such that:
\begin{enumerate}
\item \label{condtrans1} in each good box $B_\z$, one can choose a reference vertex $v_\z\in B_\z\cap V_G$,
\item \label{condtrans2} there exist $K,L$ such that to each pair of neighboring good boxes $B_{\z_1}$ and $B_{\z_2}$, one can associate a path $(v_{\z_1}, \dots,v_{\z_2})$ in $G$ between the respective reference vertices $v_{\z_1}$ and $v_{\z_2}$ of these boxes satisfying:  
\begin{enumerate}
\item $(v_{\z_1},  \dots ,v_{\z_2})\subset B_{\z_1}\cup B_{\z_2}$,
\item any edge of $(v_{\z_1}, \dots,v_{\z_2})$ has resistance at most $K$,
\item the length of $(v_{\z_1}, \dots,v_{\z_2})$ in the graph distance is bounded by $L$,
\end{enumerate}
\item \label{condtrans3} the process $X=\{X_\z, \, \z\in\ZZ^d\}:=\{\II_{B_\z \mbox{ is good}}, \, \z\in\ZZ^d\}$ stochastically dominates a supercritical independent Bernoulli site percolation process on $\ZZ^d$,
\end{enumerate}
then the random walk on $G$ with resistance $R$ is almost surely transient. 
\end{crit}
\begin{rem}\label{remtrans1}
Thanks to \cite[Theorem 0.0]{Liggett97} (see also \cite[Theorem (7.65)]{Grimmett}), in order to show \emph{(\ref{condtrans3})}, it suffices to check that $X=\{X_\z, \, \z\in\ZZ^d\}$ is a $k$-dependent process so that $\PP[X_\z=1]\geq p^*\in ]0,1[$, where $p^*$ depends on $k$ and $d$ but not on $M$. The results of Liggett, Schonmann and Stacey are still valid for processes which are not $k-$dependent but have correlations decaying rapidly. This could be used to relax the finite range of dependence assumption for the point process in Theorem \ref{proprinc}. More precisely, in order to ensure \emph{(\ref{condtrans3})}, it is enough to verify that there exists $k>0$ such that:
\[\PP \big[X_\z = 1\vert \{X_{\z'}\}_{\Vert \z'-\z\Vert\geq k}\big] \geq p\]
with $p$ as close to 1 as desired if the parameters of good boxes are well chosen. Unfortunately, we do not find any example of point process without a finite range of dependance which satisfies the assumptions of Criterion \ref{crittrans}.
\end{rem}
\begin{dem}
Thanks to condition (\ref{condtrans3}), one can define a random field $(\sigma_1,\sigma_2)\in\{0,1\}^{\ZZ^d}\times\{0,1\}^{\ZZ^d}$, where $\sigma_1$ has the distribution of the supercritical (independent) Bernoulli site percolation process, $\sigma_2$ has the law of $X$, and the pair satisfies $\sigma_1\leq\sigma_2$ almost surely. Let $\pi_\infty$ stands for the (a.s. unique) infinite percolation cluster in $\sigma_1$. It is known that the simple random walk on $\pi_\infty$ is a.s. transient when $d\geq 3$ (see \cite{Grimmett93}). By Theorem \ref{roughembbeding}, it is enough to exhibit an a.s. rough embedding from $\pi_\infty$ (with resistance 1 on each edge) to $(G,R)$. 

By stochastic domination, for any open site $\z\in\pi_\infty$ the corresponding box $B_\z$ is good. We set $\phi (\z):=v_\z$ the reference vertex of $B_\z$ given by (\ref{condtrans1}). Fix $\z_1\sim\z_2\in\pi_\infty$. Then $B_{\z_1}$ and $B_{\z_2}$ are two neighboring good boxes and one can find a path $\Phi (\z_1,\z_2)$ between $v_{\z_1}=\phi (\z_1)$ and $v_{\z_2}=\phi(\z_2)$ fully included in $B_{\z_1}\cup B_{\z_2}$. Note that $\Phi (\z_1,\z_2)$ can be assumed to be simple and that by (2)(b),(c), one has:
\[\sum_{e\in\Phi( \z_1,\z_2)}R(e)\leq \alpha :=KL.\]
Moreover, an edge of $E_G$ is in the image of at most $\beta :=2d$ edges of $\pi_\infty$ since a good box has at most $2d$ neighboring good boxes and $(v_{\z_1},\dots, v_{\z_2})\subset B_{\z_1}\cup B_{\z_2}$.
\end{dem}
\section{Recurrence in dimension 2}\label{dim2}
Let us assume that $\N$ is a stationary simple point process, almost surely in general position and satisfying {\bf (V)} and {\bf (D$_2$)}. The aim of this section is to prove the case $d=2$ of Theorem \ref{proprinc}, i.e. almost sure recurrence of walks on $\operatorname{VS}(\N)$, $\operatorname{DT}(\N)$ and $\operatorname{Gab}(\N)$. 

Since $\operatorname{Gab}(\N)$ is a subgraph of $\operatorname{DT}(\N)$, recurrence on $\operatorname{DT}(\N)$ implies recurrence on $\operatorname{Gab}(\N)$ by Rayleigh monotonicity principle. Consequently, we only need to find functions $L_{\operatorname{DT}(\N)},L_{\operatorname{VS}(\N)}$ and $N_{\operatorname{DT}(\N)},N_{\operatorname{VS}(\N)}$ so that assumptions of Criterion \ref{critrec} are satisfied for $\operatorname{DT}(\N)$ and $\operatorname{VS}(\N)$ respectively. In fact, for almost any realization of $\N$, $L(i)$ can be chosen to be of order $\sqrt{\log (i)}$ and $N(i)$ to be of order $i\sqrt{\log (i)}$. 

\subsection{Delaunay triangulation case}
We first get an upper bound on the length of the edges in $Ed_{\operatorname{DT}(\N)}(i)$ (i.e. the set of edges with only one endpoint in $B_i$) via an extended version of \cite[Lemma~1]{AddarioberrySarkar05}.

We write $\mathcal{A}_{\operatorname{DT}}(i)$ for the event `\emph{there exists an edge of $Ed_{\operatorname{DT}(\N)}(i)$ with length greater than $8c^{-\frac{1}{2}}_1\sqrt{\log i}$}' where $c_1$ is the constant appearing in {\bf (V)}, and $\mathcal{B}_{\operatorname{DT}}(i,j)$ for the event `\emph{there exists an edge of $Ed_{\operatorname{DT}(\N)}(i)$ with length between $8c^{-\frac{1}{2}}_1\sqrt{j\log i}$ and $8c^{-\frac{1}{2}}_1\sqrt{(j+1)\log i}$}'.
\begin{lemm}\label{lemmdim2}
Let $\N$ be a stationary simple point process, almost surely in general position and such that {\bf (V)} holds. Then there exists a constant $c_5>0$ such that for $i$ large enough:
\[\PP [\mathcal{A}_{\operatorname{DT}}(i)]\leq\frac{c_5}{i^2}.\]
\end{lemm} 
\begin{dem}
Assume that $\mathcal{B}_{\operatorname{DT}}(i,j)$ occurs and let $e$ be an edge of length between $8c^{-\frac{1}{2}}_1\sqrt{j\log i}$ and $8c^{-\frac{1}{2}}_1\sqrt{(j+1)\log i}$ having an endpoint in $B_i$. This edge is fully contained in $B_{i'}$, where $i':= \lceil i+8c^{-\frac{1}{2}}_1\sqrt{(j+1)\log i}\rceil$. Let $m\leq i'$ be as large as possible such that $s:=\lceil 2c^{-\frac{1}{2}}_1\sqrt{j\log i}\rceil$ divides $m$ and set $l:=m/s$. 

Now, divide $B_m$ into $l^2$ squares $Q_1,\dots ,Q_{l^2}$ of side $s$. Let $\Delta$ be one of the two Delaunay triangles having $e$ as an edge. The circumscribed sphere of $\Delta$ contains one of the two half disks $\mathfrak{D}$ with diameter $\Vert e\Vert$ located on one side or the other of $e$. Since the circumscribed sphere of the Delaunay triangle $\Delta$ contains no point of $\N$ in its interior, $\mathfrak{D}$ contains no point of $\N$ in its interior.  Note that $\mathfrak{D}$ is included in $B_{i'}$ because $e$ has an endpoint in $B_i$ and has length at most $8c^{-\frac{1}{2}}_1\sqrt{(j+1)\log i}$. Thanks to the choice of $s$, $\mathfrak{D}$ contains one of the squares $Q_1,\dots ,Q_{l^2}$, say  $Q_k$. So $Q_k\cap\N$ is empty. By stationarity of $\N$, we obtain:  

\begin{equation}
\PP [\mathcal{B}_{\operatorname{DT}}(i,j)]\leq l^2\PP[Q_1\cap\N=\emptyset ]\leq c_6i^2\PP[Q_1\cap\N=\emptyset ],\label{eqstat}
\end{equation}
for some constant $c_6>0$. But, thanks to {\bf (V)}, for $i$ large enough, one has:
\begin{align*}
\PP[Q_1\cap\N=\emptyset ]\leq\PP\left[\left[0,2c^{-\frac{1}{2}}_1\sqrt{j\log i}\right]^2\cap\N=\emptyset \right]&\leq e^{-4j\log i}=\left(\frac{1}{i^4}\right)^j.
\end{align*}
Hence, with (\ref{eqstat}):  
\begin{equation*}
\PP [\mathcal{B}_{\operatorname{DT}}(i,j)]\leq c_6i^2\left(\frac{1}{i^4}\right)^j.
\end{equation*}
Finally, for $i$ large enough, one obtains:
\begin{align*}
\PP [\mathcal{A}_{\operatorname{DT}}(i)]&\leq \sum_{j=1}^{+\infty}\PP[\mathcal{B}_{\operatorname{DT}}(i,j)]
\leq c_6i^2\sum_{j=1}^{+\infty}\left(\frac{1}{i^4}\right)^j\leq \frac{c_6i^2}{i^4-1}\leq \frac{c_5}{i^2}.
\end{align*}
\end{dem}

Thanks to the Borel-Cantelli lemma, for almost any realization of $\N$, $\mathcal{A}_{\operatorname{DT}}(i)$ holds for only finitely many $i$. Thus, for almost any realization of $\N$, one can choose $L_{\operatorname{DT}}(i)$ to be of order $\sqrt{\log i}$ in Criterion \ref{critrec}. 

We will see that we can choose $N_{\operatorname{DT}}(i)$ to be of order $i\sqrt{\log i}$. To do so, we show the following lemma and conclude with the Borel-Cantelli lemma as before.   
\begin{lemm}\label{lemmdim2b}
Let $\N$ be a stationary simple point process, almost surely in general position and such that {\bf (V)} and {\bf (D$_2$)} hold.

Then, there exists a constant $c_7>0$ such that for $i$ large enough:
\[\PP [\mathcal{C}_{\operatorname{DT}}(i)]\leq\frac{c_7}{i^2},\]
where $\mathcal{C}_{\operatorname{DT}}(i)$ is the event `\emph{$\# (Ed_{\operatorname{DT}(\N)}(i))\geq 384c^{-\frac{1}{2}}_1c_2 i\sqrt{\log i}$}'.
\end{lemm} 
\begin{dem}
Note that:
\[\PP \big[\mathcal{C}_{\operatorname{DT}}(i)\big]=\PP \big[\mathcal{C}_{\operatorname{DT}}(i)\cap \mathcal{A}_{\operatorname{DT}}(i)\big]+\PP \big[\mathcal{C}_{\operatorname{DT}}(i)\cap \mathcal{A}_{\operatorname{DT}}(i)^c\big]\leq \PP \big[\mathcal{A}_{\operatorname{DT}}(i)\big]+\PP \big[\mathcal{C}_{\operatorname{DT}}(i)\cap \mathcal{A}_{\operatorname{DT}}(i)^c\big],\]
with $\mathcal{A}_{\operatorname{DT}}(i)$ as in Lemma \ref{lemmdim2}. So, it remains to show:
\[\PP \big[\mathcal{C}_{\operatorname{DT}}(i)\cap \mathcal{A}_{\operatorname{DT}}(i)^c\big]\leq\frac{c_8}{i^2},\]
for some constant $c_8$.

On the event $\mathcal{A}_{\operatorname{DT}}(i)^c$, edges in $Ed_{\operatorname{DT}(\N)}(i)$ have lengths at most $8c^{-\frac{1}{2}}_1\sqrt{\log i}$, thus these edges are fully included in the annulus $R(i):=[-i-8c^{-\frac{1}{2}}_1\sqrt{\log i},i+8c^{-\frac{1}{2}}_1\sqrt{\log i}]^2\setminus [-i+8c^{-\frac{1}{2}}_1\sqrt{\log i},i-8c^{-\frac{1}{2}}_1\sqrt{\log i}]^2$. The restriction of $\operatorname{DT}(\N)$ to $R(i)$ is a planar graph with $\#(\N\cap R(i))$ vertices. Thanks to a corollary of Euler's formula (see \cite[Theorem 16, p.22]{ModernGraph}), it has at most $3\#(\N\cap R(i))-6$ edges. Thus $\#(Ed_{\operatorname{DT}(\N)}(i))$ is bounded from above by $3\#(\N\cap R(i))$. 

So, with {\bf (D$_2$)} and for $i$ large enough:
\begin{align*}
\PP \big[\mathcal{C}_{\operatorname{DT}}(i)\cap \mathcal{A}_{\operatorname{DT}}(i)^c\big]&\leq\PP\Big[\#\big(\N\cap R(i)\big)\geq 128c^{-\frac{1}{2}}_1c_2i\sqrt{\log i}\Big]\\
&\leq 4\PP\Big[\#\big(\N\cap \big([0,2i]\times[0,16c^{-\frac{1}{2}}_1\sqrt{\log i}]\big)\big)\geq 32c^{-\frac{1}{2}}_1c_2i\sqrt{\log i}\Big]\\
&\leq 4e^{-32c^{-\frac{1}{2}}_1c_3i\sqrt{\log i}}\leq \frac{c_8}{i^2}.
\end{align*}
\end{dem}
\subsection{Skeleton of the Voronoi tiling case}
In order to estimate the lengths of edges in $Ed_{\operatorname{VS}(\N)}(i)$, the following analogue of Lemma \ref{lemmdim2} is needed:  
\begin{lemm}\label{lemmdim2vs}
Let $\N$ be a stationary simple point process, almost surely in general position and  such that {\bf (V)} holds.

Then, there exists a positive constant $c_9$ such that for $i$ large enough:
\[\PP [\mathcal{A}_{\operatorname{VS}}(i)]\leq\frac{c_9}{i^2}\]
where $\mathcal{A}_{\operatorname{VS}}(i)$ is the event `\emph{there exists an edge of $Ed_{\operatorname{VS}(\N)}(i)$  with length greater than $4\sqrt{2}c_1^{-\frac{1}{2}}\sqrt{\log i}$}'.
\end{lemm} 
\begin{dem}
Fix $\varepsilon_0>0$ and set $\overline{B_i}:=B_i+B(0,4\sqrt{2}c_1^{-\frac{1}{2}}\sqrt{\log i}+\varepsilon_0)$ where $B(x,r)$ stands for the Euclidean ball centered at $x$ and of radius $r$. We say that $\mathcal{B}_{\operatorname{VS}}(i)$ holds if, when $\overline{B_i}$ is covered with $O(i^2/\log i)$ disjoint squares of side $2c_1^{-\frac{1}{2}}\sqrt{\log i}$ and such that one has a corner at $(-i,-i)$, each of these squares contains at least one point of $\mathcal{N}$.

We will show that, on $\mathcal{B}_{\operatorname{VS}}(i)$, every Voronoi cell intersecting $\partial B_i$ is contained in a ball centered at its nucleus and of radius $2\sqrt{2}c_1^{-\frac{1}{2}}\sqrt{\log i}$ (see Figure \ref{FigSizeVor}). This will imply that $\mathcal{B}_{\operatorname{VS}}(i)\subset \mathcal{A}_{\operatorname{VS}}(i)^c$. 

\begin{center}
\begin{figure}
\begin{center}
\includegraphics[scale=0.7]{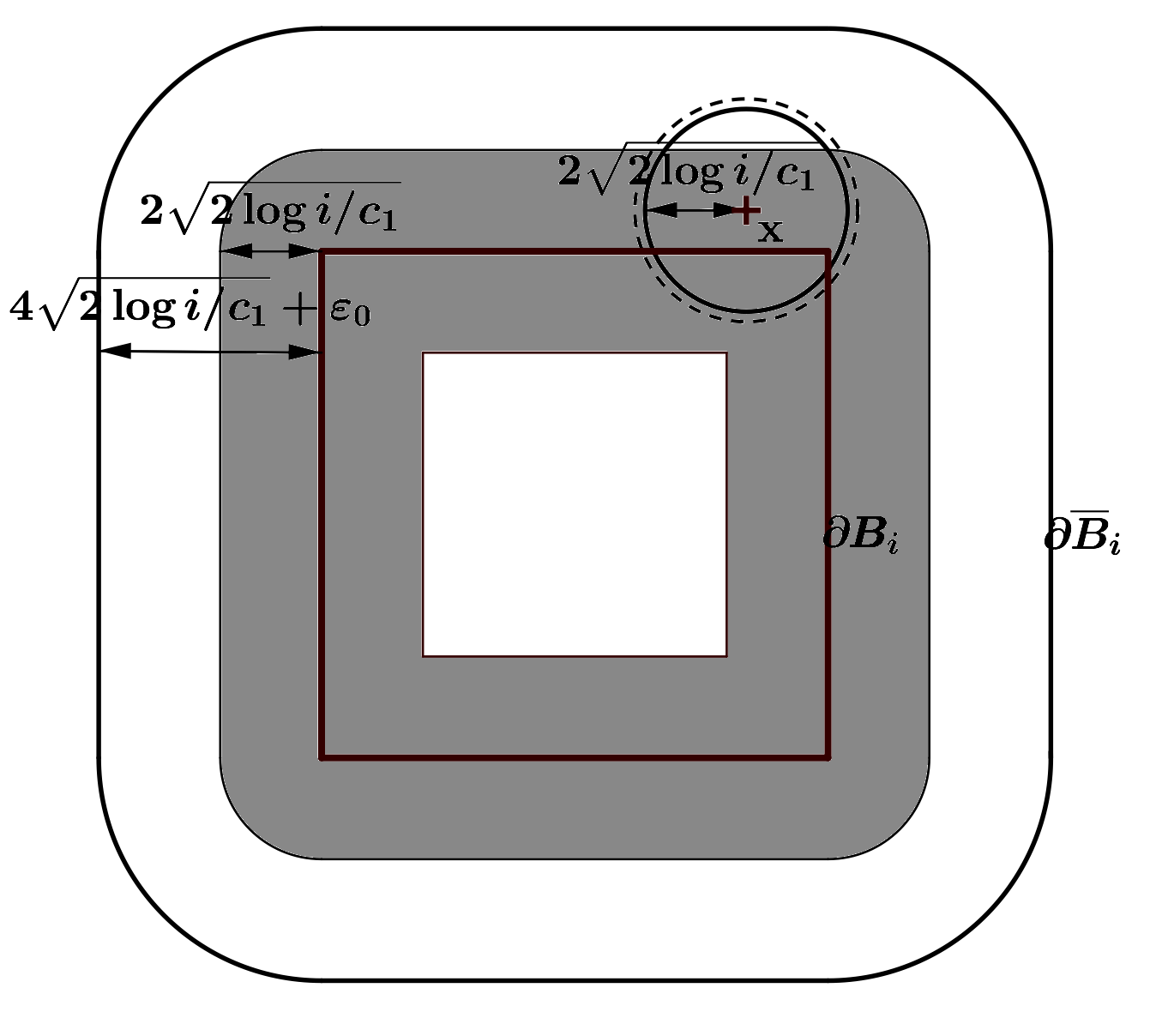}
\end{center}
\caption{\label{FigSizeVor} If $\mathcal{B}_{\operatorname{VS}}(i)$ holds, nuclei of Voronoi cells intersecting lie in the grey region. Since the dashed circle is included in $\overline{B_i}$, points of this circle are separated by at most $2\sqrt{2}c_1^{-\frac{1}{2}}\sqrt{\log i}$ from the respective nuclei of their Voronoi cells. Hence, the Voronoi cell with nucleus $\x\in\N$ is contained in $B(\x,2\sqrt{2}c_1^{-\frac{1}{2}}\sqrt{\log i})$.}
\end{figure}
\end{center}
Let us assume that $\mathcal{B}_{\operatorname{VS}}(i)$ holds. Since there is at least one point of $\N$ in each of the squares of side $2c_1^{-\frac{1}{2}}\sqrt{\log i}$ covering $\overline{B}_i$, any point in $\overline{B}_i$ has distance at most $2\sqrt{2}c_1^{-\frac{1}{2}}\sqrt{\log i}$ from the nucleus of its Voronoi cell (of each cell in which it is if it belongs to the borders of two or three Voronoi cells). In particular, nuclei of Voronoi cells intersecting $\partial B_i$ and $\partial B_i$ itself are separated by at most $2\sqrt{2}c_1^{-\frac{1}{2}}\sqrt{\log i}$. Fix $\x$ nucleus of a Voronoi cell intersecting $\partial B_i$ and note that, for $0<\varepsilon<\varepsilon_0$, $B(\x,2\sqrt{2}c_1^{-\frac{1}{2}}\sqrt{\log i}+\varepsilon)\subset\overline{B_i}$. Hence, points in $\partial B(\x,2\sqrt{2}c_1^{-\frac{1}{2}}\sqrt{\log i}+\varepsilon)$ are within a distance of at most $2\sqrt{2}c_1^{-\frac{1}{2}}\sqrt{\log i}$ from nuclei of their respective Voronoi cells which cannot be $\x$. Thus $\mathrm{Vor}_\N(\x)\subset B(\x,2\sqrt{2}c_1^{-\frac{1}{2}}\sqrt{\log i})$. So, if $\mathcal{A}_{\operatorname{VS}}(i)$ holds, $\mathcal{B}_{\operatorname{VS}}(i)$ fails. 

Finally, with {\bf (V)}, for $i$ large enough:
\begin{align*}
\PP [\mathcal{A}_{\operatorname{VS}}(i)]\leq \PP [\mathcal{B}_{\operatorname{VS}}(i)^c]&\leq c_9\frac{i^2}{\log i}\PP\big[[0,2c_1^{-\frac{1}{2}}\sqrt{\log i}]^2\cap\N=\emptyset \big]\\
&\leq c_9i^2e^{-4\log i}=\frac{c_9}{i^2}.
\end{align*}
\end{dem}

By the Borel-Cantelli lemma, for almost any realization of $\N$, $\mathcal{A}_{\operatorname{VS}}(i)$ holds only finitely many times. Thus one can choose $L_{\operatorname{VS}(\N)}$ to be of order $\sqrt{\log i}$.

It remains to show that one can choose $N_{\operatorname{VS}}(i)$ to be of order $i\sqrt{\log i}$. To do this, we state the following lemma and conclude as usual with the Borel-Cantelli lemma.    
\begin{lemm}\label{lemmdim2vsb}
Let $\N$ be a stationary simple point process, almost surely in general position and such that {\bf (V)} and {\bf (D$_2$)} hold.

Then there exists a positive constant $c_{10}$ such that for $i$ large enough:
\[\PP [\mathcal{C}_{\operatorname{VS}}(i)]\leq\frac{c_{10}}{i^2},\]
where $\mathcal{C}_{\operatorname{VS}}(i)$ is the event `\emph{$\#(Ed_{\operatorname{VS}(\N)}(i))\geq 64\sqrt{2}c_2c_1^{-\frac{1}{2}} i\sqrt{\log i}$}'.
\end{lemm} 
\begin{dem}
One has:
\[\PP \big[\mathcal{C}_{\operatorname{VS}}(i)\big]=\PP \big[\mathcal{C}_{\operatorname{VS}}(i)\cap \mathcal{B}_{\operatorname{VS}}(i)^c\big]+\PP \big[\mathcal{C}_{\operatorname{VS}}(i)\cap \mathcal{B}_{\operatorname{VS}}(i)\big]\leq \PP \big[\mathcal{B}_{\operatorname{VS}}(i)^c\big]+\PP \big[\mathcal{C}_{\operatorname{VS}}(i)\cap \mathcal{B}_{\operatorname{VS}}(i)\big],\]
with $\mathcal{B}_{\operatorname{VS}}(i)$ as in the proof of Lemma \ref{lemmdim2vs}. It remains to show that:
\[\PP \big[\mathcal{C}_{\operatorname{VS}}(i)\cap \mathcal{B}_{\operatorname{VS}}(i)\big]\leq\frac{c_{11}}{i^2},\]
for some constant $c_{11}>0$

One can see that edges in $Ed_{\operatorname{VS}(\N)}(i)$ intersect $\partial B_i$ and belong to boundaries of Voronoi cells intersecting $\partial B_i$. Note that, since Voronoi cells are convex, each of the four line segments constituting $\partial B_i$ intersects at most two sides of a given Voronoi cell and the number of edges in $Ed_{\operatorname{VS}(\N)}(i)$ intersecting one of these line segments is bounded by the number of Voronoi cells intersecting this line segment. As noticed during the proof of Lemma \ref{lemmdim2vs}, on the event $\mathcal{B}_{\operatorname{VS}}(i)$, points of $\partial B_i$ are within a distance of at most $2\sqrt{2}c_1^{-\frac{1}{2}}\sqrt{\log i}$ from nuclei of their respective Voronoi cells. Thus, if $i$ is large enough, a given Voronoi cell intersects at most two of the four line segments constituting $\partial B_i$ and nuclei of cells intersecting $\partial B_i$ are in the annulus:
\[R(i):=\Big[-i-2\sqrt{2}c_1^{-\frac{1}{2}}\sqrt{\log i},i+2\sqrt{2}c_1^{-\frac{1}{2}}\sqrt{\log i}\Big]^2\setminus \Big[-i+2\sqrt{2}c_1^{-\frac{1}{2}}\sqrt{\log i},i-2\sqrt{2}c_1^{-\frac{1}{2}}\sqrt{\log i}\Big]^2.\] So, on $\mathcal{B}_{\operatorname{VS}}(i)$, for $i$ large enough, the number of edges in $Ed_{\operatorname{VS}(\N)}(i)$ is bounded by twice the number of Voronoi cells intersecting $\partial B_i$ which is less than $2\#(\N\cap R(i))$.

Thus, with {\bf (D$_2$)}, for $i$ large enough:
\begin{align*}
\PP \big[\mathcal{C}_{\operatorname{VS}}(i)\cap \mathcal{B}_{\operatorname{VS}}(i)\big]&\leq\PP\big[\#\big(\N\cap R(i)\big)\geq 32\sqrt{2}c_2c_1^{-\frac{1}{2}}i\sqrt{\log i}\big]\\
&\leq 4\PP\Big[\#\big(\N\cap \big([0,2i]\times [0,4\sqrt{2}c_1^{-\frac{1}{2}}\sqrt{\log i}]\big)\big)\geq 8\sqrt{2}c_2c_1^{-\frac{1}{2}}i\sqrt{\log i}\Big]\\
&\leq 4 e^{-8\sqrt{2}c_3c_1^{-\frac{1}{2}}i\sqrt{\log i}}\leq \frac{c_{11}}{i^2}.
\end{align*}
\end{dem}

\section{Transience in higher dimensions}\label{dim3}
This section is devoted to the proof of the second part of Theorem \ref{proprinc}, i.e. almost sure transience of walks on $\operatorname{VS}(\N)$, $\operatorname{DT}(\N)$ and $\operatorname{Gab}(\N)$ under the assumptions on $\N$.

The assumption that $\N$ has almost surely no descending chain is only used in the proof of transience on $\operatorname{Gab}(\N)$. Since $\operatorname{Gab}(\N)$ is a subgraph of $\operatorname{DT}(\N)$, transience on $\operatorname{Gab}(\N)$ implies transience on $\operatorname{DT}(\N)$. Actually, one can directly prove transience on $\operatorname{DT}(\N)$ without this additional condition. This proof is very similar to the $\operatorname{VS}(\N)$ case and is omitted here. 

We will define `good boxes' and use Criterion \ref{crittrans} for $\operatorname{Gab}(\N)$ and $\operatorname{VS}(\N)$: we construct paths needed in Criterion \ref{crittrans} and check that $\PP[X_\z=1]$ is large enough if parameters are well chosen. Thanks to the finite range of dependence assumptions (say of range $k$) and to the definitions of `good boxes' it will be clear that $\{X_\z\}$ is a $k$-dependent Bernoulli process on $\ZZ^d$.
\subsection{Skeleton of the Voronoi tiling case}
\subsubsection{Good boxes}\label{GBVor}
For $M\geq 1$ to be determined later, consider a partition of $\RR^d$ into boxes of side $M$:
\[B_\z=B^M_\z:=M\z+\Big[-\frac{M}{2},\frac{M}{2}\Big)^d,\,\z\in \ZZ^d.\]
We say that a box $B_\z$ is $M-$good for $\operatorname{VS}(\N)$ if the following conditions are satisfied:
\begin{itemize}
\item[{\it -i-}] $\# \big(B_\z\cap\N\big)\leq 2c_4M^d$,
\item[{\it -ii-}] when $B_\z$ is (regularly) cut into $\alpha_d^d:=(6\lceil\sqrt{d}\rceil)^d$ sub-boxes $b_i^\z$ of side $M/\alpha_d$, each of these sub-boxes contains at least one point of $\N$.
\end{itemize}
 
\subsubsection{Construction of paths}
First of all, for each good box $B_{\z_i}$, one can fix a reference vertex $v_{\z_i}$ on the boundary of the Voronoi cell of $M\z_i$ as required by condition (\ref{condtrans1}) in Criterion~\ref{crittrans}. Note that, thanks to {\it -ii-} in the definition of good boxes, the cell containing $M\z_i$ is included in $B_{\z_i}$ so that $v_{\z_i}\in B_{\z_i}$.

Consider two neighboring good boxes $B_{\z_1}$ and $B_{\z_2}$. One must show that there exists a path in $\operatorname{VS}(\N)$ between $v_{\z_1}$ and $v_{\z_2}$ so that condition (\ref{condtrans2}) of Criterion \ref{crittrans} holds. To do so, note that one can find a self-avoiding path between $v_{\z_1}$ and $v_{\z_2}$ with edges belonging to boundaries of the Voronoi cells crossing the line segment $[M\z_1,M\z_2]$. Clearly, such a path is contained in $B_{\z_1}\cup B_{\z_2}$ as soon as Voronoi cells intersecting $[M\z_1,M\z_2]$ are included in $B_{\z_1}\cup B_{\z_2}$. Since sub-boxes $b^{\z_1}_1,\dots ,b^{\z_1}_{\alpha_d^d},b^{\z_2}_1,\dots ,b^{\z_2}_{\alpha_d^d}$ are non-empty (of points of $\N$), points in $B_{\z_1}\cup B_{\z_2}$ are distant by at most $\sqrt{d}M/\alpha_d$ from nuclei of their respective Voronoi cell. As in the proof of Lemma \ref{lemmdim2vs}, one can see that the Voronoi cells intersecting $[M\z_1,M\z_2]$ are included in $[M\z_1,M\z_2]+B(0,2\sqrt{d}M/\alpha_d)\subset B_{\z_1}\cup B_{\z_2}$. Actually, since there is at least one point of $\N$ in each of the sub-boxes of side  $M/\alpha_d$ contained in $B_{\z_1}\cup B_{\z_2}$, nuclei of Voronoi cells intersecting $[M\z_1,M\z_2]$ are within a distance of at most $\sqrt{d}M/\alpha_d$ from this line segment. Fix $\x$  the nucleus of a Voronoi cell intersecting $[M\z_1,M\z_2]$  and note that, for $\varepsilon>0$ small enough, $B(\x,\sqrt{d}M/\alpha_d+\varepsilon)\subset B_{\z_1}\cup B_{\z_2}$. Hence, points in $\partial B(\x,\sqrt{d}M/\alpha_d+\varepsilon)$ are within a distance of at most $\sqrt{d}M/\alpha_d$ from the nuclei of their respective Voronoi cells which cannot be $\x$. Repeating above arguments with $[M\z_1,M\z_2]+B(0,2\sqrt{d}M/\alpha_d)$ instead of $[M\z_1,M\z_2]$, one obtains that Voronoi cells intersecting $[M\z_1,M\z_2]+B(0,2\sqrt{d}M/\alpha_d)$ have their nuclei in $[M\z_1,M\z_2]+B(0,3\sqrt{d}M/\alpha_d)$. Since Voronoi cells intersecting $[M\z_1,M\z_2]$ are contained in $[M\z_1,M\z_2]+B(0,2\sqrt{d}M/\alpha_d)$, it follows that nuclei of Voronoi cells which are neighbors of cells intersecting $[M\z_1,M\z_2]$ are in $[M\z_1,M\z_2]+B(0,3\sqrt{d}M/\alpha_d)\subset B_{\z_1}\cup B_{\z_2}$. Note that the path between $v_{\z_1}$ and $v_{\z_2}$ has chemical length (i.e. for the graph distance) bounded by the number of vertices on the boundaries of Voronoi cells crossing $[M\z_1,M\z_2]$. This is less than the number of Voronoi cells crossing $[M\z_1,M\z_2]$ times the maximal number of vertices on such a cell. Since these cells are included in $B_{\z_1}\cup B_{\z_2}$ and $B_{\z_1},\,B_{\z_2}$ are good, there are no more than $4c_4M^d$ cells intersecting $[M\z_1,M\z_2]$. Each of these cells has at most $4c_4M^d$ neighboring cells whose nuclei are in $B_{\z_1}\cup B_{\z_2}$. The total number of vertices on the boundary of such a cell is generously bounded by $\binom{4c_4M^d}{d}$. Indeed, any of these vertices is obtained as the intersection of $d$ bisecting hyperplanes separating the cell from one of its neighbors (which are at most $4c_4M^d$). So, one can choose $L:=4c_4M^d\binom{4c_4M^d}{d}$ in Criterion \ref{crittrans} (2)(c). 

Finally, if $C$ is uniformly bounded from below set $K:=\max 1/C$ in Criterion \ref{crittrans} (2)(b). If $C$ is given by a decreasing positive function $\varphi$ of edge length, set $K:=1/\varphi(\sqrt{d+3}M)$. Indeed, by construction, any edge in the path from $v_{\z_1}$ to $v_{\z_1}$ has (Euclidean) length at most $\sqrt{d+3}M$ since it is included in $B_{\z_1}\cup B_{\z_2}$. 
\subsubsection{$\PP[X_\z=1]$ is large enough}
It remains to show that if $M$ is appropriately chosen (3) in Criterion \ref{crittrans} is satisfied. Since $\N$ has a finite range of dependence $k$, assuming that $M\geq 2$, the process $\{X_\z\}$ is $k-$dependent, so, as noticed in Remark \ref{remtrans1}, it suffices to show that for $M$ large enough:
\[\PP\big[X_\z=1\big]\geq p^*\]
where $p^*=p^*(d,k)<1$ is large enough to ensure that $\{X_\z\}$ dominates site percolation on $\ZZ^d$.

Indeed, with {\bf (V)} and {\bf (D$_{3^+}$)}, one has:
\begin{align*}
\PP\big[X_\z=0\big]&\leq \PP\Big[\# \big(B_\z\cap\N\big)>2c_4M^d\mbox{ or one of the }b^\z_i\mbox{s is empty}\Big]\\
&\leq \PP\Big[\# \big([0,M]^d\cap\N\big)>2c_4M^d\Big]+\alpha_d^d\PP\Big[\# \Big(\Big[0,\frac{M}{\alpha_d}\Big]^d\cap\N\Big)=0\Big]\\
&\leq \exp \big(-c_4M^d\big)+\big(6\lceil\sqrt{d}\rceil\big)^d\exp \Big(-c_1\frac{M^d}{6^d\lceil\sqrt{d}\rceil^d}\Big) 
\end{align*}
which is as small as we wish for large $M$.
\subsection{Gabriel graph case}
\subsubsection{A geometric lemma} We shall state a generalization of \cite[Lemma 1]{BBD} which allows us to control the behavior of paths as needed in Criterion \ref{crittrans}. 
\begin{lemm}\label{lemmgeo}Let $\N$ be a locally finite subset of $\RR^d$ without descending chains.

Then, for any $x,y\in \N$, there exists a path $(x_1=x, \dots,x_n=y)$ in $\operatorname{Gab}(\N)$ from $x$ to $y$ such that:
\begin{equation}
\sum_{i=1}^{n-1}\Vert x_{i+1}-x_i\Vert^2\leq \Vert y-x\Vert^2 .\label{Gabpath}
\end{equation}  
\end{lemm}
\begin{dem}
Assume that $\N$ is a locally finite subset of $\RR^d$ in which there are $x,y\in\N$ such that there is no path between $x$ and $y$ in $\operatorname{Gab}(\N)$ satisfying (\ref{Gabpath}). We fix such $x$ and $y$ and we obtain a contradiction by constructing a descending chain $(u_i)_{i\in\NN}$ in $\N$. To this end, we prove by induction that for any $i\geq 1$ there exist $u_0,\dots ,u_i, z_1^{i+1}\in\N$ such that:
\begin{enumerate}[-i-]
\item \label{RECGAB1}$\Vert u_{j+1}-u_j\Vert<\Vert u_j-u_{j-1}\Vert$ for $j=1,\dots,i-1$,
\item $\Vert z_1^{i+1}-u_i\Vert< \Vert u_i-u_{i-1}\Vert$,
\item \label{RECGAB3}there is no path between $u_i$ and $z^{i+1}_1$ satisfying (\ref{Gabpath}).
\end{enumerate}

This eventually proves that $\N$ has a descending chain $(u_i)_{i\in\NN}$.
\paragraph{\emph{Base step.}} Since there is no path between $x$ and $y$ in $\operatorname{Gab}(\N)$ satisfying (\ref{Gabpath}), the open ball $B([x,y])$ with diameter $[x,y]$ contains at least a point of $\N$. Note that, otherwise, there would be an edge between $x$ and $y$ which would be a particular path between these points trivially satisfying (\ref{Gabpath}). Let us denote by $z^1_2$ the smallest point in the lexicographic order in $B([x,y])\cap\N$. Then, either there is no path satisfying (\ref{Gabpath}) between $x$ and $z^1_2$ or between $z^1_2$ and $y$. Indeed, the section of $B([x,y])$ by the plane generated by $x,y,z^1_2$ is a circle with diameter $[x,y]$. Thus we have:
\begin{equation*}
\Vert z^1_2-x\Vert^2+\Vert y-z^1_2\Vert^2\leq\Vert y-x\Vert^2.\label{Alkashi}
\end{equation*}
If there were a path $(x_1=x,\dots ,x_k=z^1_2)$ between $x$ and $z^1_2$ and a path $(x_k=z,\dots ,x_n=y)$ between $z^1_2$ and $y$ satisfying (\ref{Gabpath}), then their concatenation would be a path between $x$ and $y$ satisfying (\ref{Gabpath}):
\begin{align*}
\sum_{i=1}^{n-1}\Vert x_{i+1}-x_i\Vert^2&=\sum_{i=1}^{k-1}\Vert x_{i+1}-x_i\Vert^2+\sum_{i=k}^{n-1}\Vert x_{i+1}-x_i\Vert^2\\
&\leq \Vert z^1_2-x \Vert^2+\Vert y-z^1_2 \Vert^2\leq\Vert y-x \Vert^2.
\end{align*}
Note also that $\Vert z^1_2-x \Vert ,\Vert y-z^1_2 \Vert<\Vert y-x \Vert$. 
Without loss of generality, we may assume that there is no path satisfying (\ref{Gabpath}) between $z^1_2$ and $y$ and we set $u_0:=x$, $u_1:=y$ and $z^2_1:=z^1_2$.
       
\paragraph{\emph{Inductive step.}} We assume that $u_0,\dots, u_i, z_1^{i+1}\in\N$ satisfying properties -\ref{RECGAB1}- to -\ref{RECGAB3}- listed above are constructed, and we construct $u_{i+1}$ and $z_1^{i+2}$ such that:
\begin{enumerate}[-i-]
\item $\Vert z_1^{i+2}-u_{i+1}\Vert<\Vert u_{i+1}-u_i\Vert<\Vert u_i-u_{i-1}\Vert$,
\item there is no path between $u_{i+1}$ and $z^{i+2}_1$ satisfying (\ref{Gabpath}).
\end{enumerate}
Since there is no path between $u_i$ and $z^{i+1}_1$ in $\operatorname{Gab}(\N)$ satisfying (\ref{Gabpath}), $B([u_i,z^{i+1}_1])\cap \N\neq\emptyset$ (otherwise, there was an edge between $u_i$ and $z^{i+1}_1$ in $\operatorname{Gab}(\N)$). We write $z^{i+1}_2$ for the smallest point in the lexicographic order in $B([u_i,z^{i+1}_1])\cap\N$. There are then two possibilities:
\begin{description}
\item[1) there is no path between $z^{i+1}_1$ and  $z^{i+1}_2$ satisfying (\ref{Gabpath})] We then set $u_{i+1}:=z^{i+1}_1$ and $z_1^{i+2}:=z_2^{i+1}$. Using that $z^{i+2}_1$ is inside the ball $B([u_i,u_{i+1}])$ and the induction hypothesis, we have:
\[\Vert z_1^{i+2}-u_{i+1}\Vert <\Vert u_{i+1}-u_i\Vert =\Vert z_1^{i+1}-u_i\Vert <\Vert u_i-u_{i-1}\Vert <\dots <\Vert u_1-u_0\Vert.\] 
Moreover, there is no path between $u_{i+1}$ and $z^{i+2}_1$ satisfying (\ref{Gabpath}).
\item[2)  there is a path between $z^{i+1}_1$ and  $z^{i+1}_2$ satisfying (\ref{Gabpath})] Then there is no path between $u_i$ and $z^{i+1}_2$  satisfying (\ref{Gabpath}), otherwise there would be a path between $u_{i}$ and $z^{i+1}_1$ satisfying (\ref{Gabpath}). In particular, $B([u_i,z^{i+1}_2])\cap\N\neq\emptyset$ and we denote by $z^{i+1}_3$ the smallest point in the lexicographic order in $B([u_i,z^{i+1}_2])\cap\N$. If there is no path between $z^{i+1}_3$ and $ z^{i+1}_2$ satisfying (\ref{Gabpath}), we proceed as in {\bf 1)}. Otherwise, we repeat the procedure until we find $z^{i+1}_{n_{i+1}}$ and $z^{i+1}_{n_{i+1}+1}$ such that there is no path between $z^{i+1}_{n_{i+1}}$ and $z^{i+1}_{n_{i+1}+1}$ in $\operatorname{Gab}(\N)$ satisfying (\ref{Gabpath}). One can see that:
\[\Vert z^{i+1}_{n_{i+1}+1}-z^{i+1}_{n_{i+1}}\Vert<\Vert z^{i+1}_{n_{i+1}}-u_i\Vert<\dots<\Vert z_1^{i+1}-u_i\Vert<\Vert u_i-u_{i-1}\Vert.\]
In particular, such $n_{i+1}$ does exist. Indeed, if not, there would be infinitely many points of $\N$ in the ball $B(u_i,\Vert u_i-u_{i-1}\Vert)$ and $\N$ would not be locally finite. We set $u_{i+1}:=z^{i+1}_{n_{i+1}}$ and $z_1^{i+2}:=z_{n_{i+1}+1}^{i+1}$. Then:
\[\Vert z_1^{i+2}-u_{i+1}\Vert <\Vert u_{i+1}-u_i\Vert <\Vert u_i-u_{i-1}\Vert <\dots <\Vert u_1-u_0\Vert,\] 
and there is no path between $u_{i+1}$ and $z^{i+2}_1$ satisfying (\ref{Gabpath}). 
\end{description}
Finally, $\N$ has a descending chain. This contradicts the original assumption.
\end{dem}
\subsubsection{Good boxes}\label{GBGab}
For $M\geq 1$, consider as before a partition of $\RR^d$ into boxes of side $M$, $\{B_\z,\,\z\in\ZZ^d\}$. For $m\in\NN^*$, write $\alpha_{d,m}$ for the odd integer such that:
\[\beta_dm^2+\sqrt{d}+1\leq\alpha_{d,m}<\beta_dm^2+\sqrt{d}+3,\]
where $\beta_d:=2^{2d+2}\big(d+3+2(d+3)^{3/2}\big)$. We say that a box $B_\z$ is $(M,m)-$good for the Gabriel graph if when $B_\z$ is cut into $\alpha_{d,m}^d$ sub-boxes $b_i^\z$ of side $M/\alpha_{d,m}$, each of these sub-boxes contains at least one and at most $m$ points of $\N$.
\subsubsection{Construction of paths}
Let $B_{\z_1}$ and $B_{\z_2}$ be two neighboring good boxes of side $M$ and $b_1^{\z_1,\z_2},\dots,b_{\alpha_{d,m}+1}^{\z_1,\z_2}$ be the sub-boxes of side $M/\alpha_{d,m}$ intersecting the line segment $[M\z_1,M\z_2]$. Write $c_i$ (resp. $v_i$) for the center of $b_i^{\z_1,\z_2}$ (resp. the point of $b_i^{\z_1,\z_2}\cap\N$ which is the closest to $c_i$). Vertices $v_1$ and $v_{\alpha_{d,m}+1}$ are reference vertices of $B_{\z_1}$ and $B_{\z_2}$ respectively. One must prove that there exists a Gabriel path (i.e. a path in the Gabriel
graph) from $v_1$ to $v_{\alpha_{d,m}+1}$ which is included in $B_{\z_1}\cup B_{\z_2}$. To this end, it suffices to check that, for all $i\in\{1,\dots, \alpha_{d,m}\}$, there exists a Gabriel path from $v_i$ to $v_{i+1}$ which is included in $B(c_i,M/2)$. Proceeding along the same lines as in the proof of \cite[Lemma~3]{BBD}, we show that the path $\gamma_i:=(x_1=v_i, \dots, x_n=v_{i+1})$ given by Lemma \ref{lemmgeo} for $v_i$ and $v_{i+1}$ satisfies this property. Observe that this path can be assumed to be simple and satisfies: 
\begin{equation}
\sum_{j=1}^{n-1}\Vert x_{j+1}-x_j\Vert^2\leq\Vert v_{i+1}-v_i\Vert^2\leq (d+3)\frac{M^2}{\alpha_{d,m}^2},\label{eqchem1}
\end{equation}
where the last bound is obtained by the Pythagorean theorem using that $v_i$ and $v_{i+1}$ belong to two neighboring boxes of side $M/\alpha_{d,m}$. 
In particular, $\gamma_i$ does not contain any edge with length greater than $\sqrt{d+3}M/\alpha_{d,m}$ and contains at most $2^{2d+2}(d+3)m^2$ edges with length between $M/(2^{d+1}\alpha_{d,m}m)$ and $\sqrt{d+3}M/\alpha_{d,m}$, called \emph{long} edges in the following.   
Indeed, with (\ref{eqchem1}), 
\begin{align*}
\frac{M^2}{2^{2d+2}\alpha^2_{d,m}m^2}\#\big\{e:\mbox{ long edge of } \gamma_i\big\}
&\leq\sum_{e:\mbox{ \tiny long edge}}\Vert e\Vert^2
\leq (d+3)\frac{M^2}{\alpha_{d,m}^2}.
\end{align*}

Hence, $\gamma_i$ consists of at most $2^{2d+2}(d+3)m^2$ long edges and at most $2^{2d+2}(d+3)m^2+\nolinebreak1$ groups of consecutive \emph{short} edges with length lower than $M/(2^{d+1}\alpha_{d,m}m)$). Let us consider a simple path starting at $v_i$ and having $N_l\in\{0,\dots,2^{2d+2}(d+3)m^2\}$ long edges and $N_s\in\{0,\dots,2^{2d+2}(d+3)m^2+1\}$ groups of consecutive short edges. Each group of short egdes has total Euclidean length bounded by $M/(2\alpha_{d,m})$. Indeed, one can prove by induction that the vertex at the begining of each group of consecutive short edges satisfies the assumption of Lemma~\ref{fait}, stated at the end of this subsection for sake of readability. Consequently, the distance from $c_i$ to the farthest point of the path is less than:
\begin{align*}
\Vert v_i -c_i\Vert &+N_l\frac{\sqrt{d+3}M}{\alpha_{d,m}}+N_s\frac{M}{2\alpha_{d,m}}\\
&\leq \sqrt{d}\frac{M}{2\alpha_{d,m}} +N_l\frac{\sqrt{d+3}M}{\alpha_{d,m}}+N_s\frac{M}{2\alpha_{d,m}}\\
&\leq \sqrt{d}\frac{M}{2\alpha_{d,m}} + 2^{2d+2}(d+3)^\frac{3}{2}m^2\frac{M}{\alpha_{d,m}}
+ \big(2^{2d+2}(d+3)m^2+1\big)\frac{M}{2\alpha_{d,m}}\\
&=\frac{M}{2}\frac{\beta_d m^2+\sqrt{d}+1}{\alpha_{d,m}}\leq\frac{M}{2}.
\end{align*}      
Thus $\gamma_i$ is included in $B(c_i,M/2)$ and there exists a Gabriel path $\gamma:=(v_1,\dots, v_{\alpha_{d,m}+1})$ from $v_1$ to $v_{\alpha_{d,m}+1}$ contained in $B_{\z_1}\cup B_{\z_2}$.

It remains to choose $L$ and $K$ in Criterion \ref{crittrans} (2)(b) and (2)(c). Since $\gamma\subset B_{\z_1}\cup B_{\z_2}$ and can be supposed simple, it has chemical length at most $\#\big((B_{\z_1}\cup B_{\z_2})\cap\N\big)-1\leq 2\alpha_{d,m}^dm-1$. Thus, one can set $L:=2\alpha_{d,m}^dm-1$. If $C(\cdot)=\varphi (\Vert\cdot\Vert)$ is a decreasing positive function of edges lengths, we set $K:=1/\varphi(\sqrt{d+3}M)$; if $C$ is uniformly bounded from below, we set $K:=\max 1/C$.

\begin{lemm}\label{fait}
Let $B_{\z_1},\, B_{\z_2}$ be two neighboring good boxes and  $\gamma$ be a simple Gabriel path with edges of length bounded by $M/(2^{d+1}\alpha_{d,m}m)$ with a vertex $u$ such that
$B(u,M/(2\alpha_{d,m}))\subset B_{\z_1}\cup B_{\z_2}$.

Then, $\gamma$ consists in at most $2^dm-1$ edges. In particular, it has total Euclidean length bounded by $M/(2\alpha_{d,m})$.
\end{lemm}  
\begin{dem}
Note that $B(u,M/(2\alpha_{d,m}))$ intersects at most $2^d$ sub-boxes of side $M/\alpha_{d,m}$. Since it is furthermore included in $B_{\z_1}\cup B_{\z_2}$ which are good boxes, it contains at most $2^dm$ points of $\N$. Assume that $\gamma$ contains more than $2^dm$ edges. Hence, $\gamma$ has a (sub-)path $\gamma'$ of $2^dm$ edges (and $2^dm+1$ vertices) such that $u\in\gamma'$. But $\gamma'$ is included in $B(u,M/(2\alpha_{d,m}))$ which provides a contradiction. 
\end{dem}
\subsubsection{$\PP[X_\z=1]$ is large enough}
Assuming that $M\geq 2$, the process $\{X_\z\}$ is $k-$dependent because of the definition of good boxes and the fact that $\N$ has a finite range of dependence. It remains to show that if $M,m$ are suitably chosen:
\[\PP\big[ X_\z=1\big]\geq p^*\]
where $p^*=p^*(d,k)<1$ is large enough to ensure that $\{X_\z\}$ dominates supercritical site percolation on $\ZZ^d$.

Thanks to the choice of $\alpha_{d,m}\sim \beta_dm^2$, for $m\in\NN^*$ large enough, we can choose $M$ so that:
\begin{equation*}
\frac{\alpha_{d,m}^d}{c_1}\log\Bigg(\frac{2\alpha_{d,m}^d}{1-p^*}\Bigg)\leq M^d\leq\frac{\alpha_{d,m}^d}{c_4} \Bigg(m-\log\Bigg(\frac{2\alpha_{d,m}^d}{1-p^*}\Bigg)\Bigg).\label{eqM}
\end{equation*}

With {\bf (V)} and {\bf (D$_{3^+}$)}:
\[\PP\Big[\#\Big(\Big[0,\frac{M}{\alpha_{d,m}}\Big]^d\cap\N\Big)=0\Big]\leq\exp\Big(-c_1\frac{M^d}{\alpha_{d,m}^d}\Big)\leq\frac{1-p^*}{2\alpha_{d,m}^d},\]
and
\[\PP\Big[\#\Big(\Big[0,\frac{M}{\alpha_{d,m}}\Big]^d\cap\N\Big)>m\Big]\leq\exp\Big(c_4\frac{M^d}{\alpha_{d,m}^d}-m\Big)\leq\frac{1-p^*}{2\alpha_{d,m}^d}.\]

Finally, using stationarity of $\N$:
\begin{align*}
\PP[X_\z=0]&\leq\sum_{i=1}^{\alpha_{d,m}^d}\Big\{\PP\big[\#(b^\z_i\cap\N)=0\big]+\PP\big[\#(b^\z_i\cap\N)>m\big]\Big\}\\
&=\alpha_{d,m}^d\Big\{\PP\Big[\#\Big(\Big[0,\frac{M}{\alpha_{d,m}}\Big]^d\cap\N\Big)=0\Big]+\PP\Big[\#\Big(\Big[0,\frac{M}{\alpha_{d,m}}\Big]^d\cap\N\Big)>m\Big]\Big\}\\
&\leq 1-p^*.
\end{align*}
\section{Examples of point processes}\label{exproc}
In this section, an overview on point processes appearing in Theorem \ref{thprinc} is given. Assumptions of Theorem \ref{proprinc} are checked for these processes. Note that the probability estimates given here are rough but good enough to verify {\bf (V)}-{\bf (D$_{3^+}$)}. 
\subsection{Poisson point processes}
For homogeneous Poisson point processes (PPPs), stationarity and the finite range of dependence condition are clear. Note that stationary Poisson point processes are almost surely in general position. We refer to \cite{HM96} or \cite{DL05} for the almost sure absence of descending chains. Moreover, we check by standard computations that, for any Borel set $A$, a PPP $\N$ of intensity $\lambda$ satisfies:
\begin{equation}\label{EqPPPVoid}
\PP\big[\#\big(A\cap\N\big)=0\big]=e^{-\lambda \operatorname{Vol}_{\RR^d}(A)},
\end{equation}
and
\begin{equation}\label{EqPPPDev}
\PP\big[\#\big(A\cap\N\big)>m\big]\leq e^{\lambda(e-1) \operatorname{Vol}_{\RR^d}(A)-m}.
\end{equation}
This implies {\bf (V)}-{\bf (D$_{3^+}$)}.

Results similar to (\ref{EqPPPVoid}) and (\ref{EqPPPDev}) are satisfied when the intensity measure $\mu$ of the PPP is comparable to Lebesgue measure on $\RR^d$ in the sense that there exists a positive constant $c_{12}$ such that for every measurable subset $A$ of $\RR^d$:
\[\frac{1}{c_{12}} \operatorname{Vol}_{\RR^d}(A)\leq \mu (A) \leq c_{12}\operatorname{Vol}_{\RR^d}(A).\]
This implies that the conclusions of Theorem \ref{thprinc}  also hold for such non-stationary PPPs (see Remark \ref{remNonStat}).
\subsection{Mat\'ern cluster processes}
Mat\'ern cluster processes (MCPs) are particular cases of Neyman-Scott Poisson processes (see \cite[p. 171]{ChiuSKM}). Cluster processes are used as models for spatial phenomena, e.g. galaxy locations in space \cite{Galaxy} or epicenters of micro-earthquake locations \cite{VJEarthquakes}. 

MCPs are constructed as follows. Let $\lambda, \mu, R >0$. One first chooses a PPP $Y$ of intensity $\lambda$ called the \emph{parent process}. For any $y\in Y$, a centered \emph{daughter process} $\N_y$ is then chosen such that, given $Y$, $\{\N_y\}_{y\in Y}$ are mutually independent PPP  with intensity $\mu$ in $B(0,R)$. Then, $\N:=\bigcup_{y\in Y}(y+\N_y)$ is a MCP with parameters $\lambda, \mu, R$. It is clear that such processes are stationary. Since its parent process has a finite range of dependence and daughter processes have bounded supports, any MCP has a finite range of dependence. Thanks to \cite[Proposition 2.3]{HNS12} or \cite[Theorem 7.2]{DL05}, MCPs have almost surely no descending chains.
MCPs can be seen as doubly stochastic processes or Cox processes (see \cite[p. 166]{ChiuSKM} for a definition). Their (diffusive) random intensity measures $\mu_Y$ have densities $\sum_{y\in Y}\II_{B(y,R)}(x)$ w.r.t. Lebesgue measure. In particular, $(d-1)$ dimensional hyperplanes and spheres are $\mu_Y-$null sets and MCPs are almost surely in general position.
   
It remains to check assumptions {\bf (V)}-{\bf (D$_{3^+}$)}. Let $L>2R$, then,
\begin{align*}
\PP\big[\#\big(\N\cap [0,L]^d \big)=0\big]&=\EE\big[\PP\big[\#\big(\N\cap [0,L]^d \big)=0\big\vert Y\big]\big]\\
&=\EE\Big[\PP\Big[\underset{y\in Y}{\bigcap}\big\{\#\big(\N_y\cap [0,L]^d \big)=0\big\}\Big\vert Y\Big]\Big]\\
&=\EE\Big[\prod_{y\in Y}\PP\Big[\big\{\#\big(\N_y\cap [0,L]^d \big)=0\big\}\Big\vert Y\Big]\Big]\\
&=\EE\Big[\prod_{y\in Y}\exp\big(-\mu\operatorname{Vol}_{\RR^d}(B(y,R)\cap [0,L]^d)\big)\Big]\\
\end{align*}
so, with \cite[Theorem 3.2.4]{SW}, 
\begin{align*}
\PP\big[\#\big(\N\cap [0,L]^d\big)=0\big]
&=\exp\Big(-\lambda\int_{[0,L]^d+B(0,R)}\hspace{-1cm}\big\{1-\exp(-\mu\operatorname{Vol}_{\RR^d}([0,L]^d\cap B(y,R)))\big\}\d y\Big)\\
&\leq \exp\Big(-\lambda\int_{[0,L]^d-B(0,R)}\hspace{-1cm}\big\{1-\exp(-\mu\operatorname{Vol}_{\RR^d}(B(0,R)))\big\}\d y\Big)\\
\end{align*}
where, for $A, B\subset \RR^d$,  $A-B:=\{x\in \RR^d :\,\forall y\in B, \, x+y\in A\} $. Thus,
\begin{align*}
\PP\big[\#\big(\N\cap [0,L]^d \big)=0\big]
&\leq \exp\Big(-\lambda\big\{1-\exp(-\mu\operatorname{Vol}_{\RR^d}(B(0,R)))\big\}\big(L-2R\big)^d\Big)
\end{align*}
and {\bf (V)} holds. Alternatively, one can see that the MCP with parameters $\lambda,\mu$ and $R$ satisfies the hypothesis
 {\bf (V)} because it stochastically dominates a PPP with intensity $\lambda'$ where $\lambda' =\lambda (1-\exp(\mu \operatorname{Vol}_{\RR^d}(B(0,R)))$, obtained by deleting all but one uniformly chosen point of each non-empty daughter process $\mathcal{N}_y$.

In order to show {\bf (D$_2$)} and {\bf (D$_{3^+}$)} for MCPs one can use exponential Markov inequality and the following estimate. For any bounded Borel set $A$:
\begin{align*}
\EE\big[e^{\# (\N\cap A)}\big]
&=\EE\Big[\EE\Big[\exp\Big(\sum_{y\in Y}\# (\N_y\cap A)\Big)\Big\vert Y\Big]\Big]\\
&=\EE\Big[\prod_{y\in Y\cap (A+B(0,R))}\EE\Big[e^{\# (\N_y\cap A)}\Big\vert Y\Big]\Big]\\
&\leq\EE\Big[\prod_{y\in Y\cap (A+B(0,R))}\EE\Big[e^{\#\N_y}\Big\vert Y\Big]\Big]\\
&=\EE\Big[\exp\Big(\mu\operatorname{Vol}_{\RR^d}(B(0,R))(e-1)\# \big\{y\in Y\cap (A+B(0,R))\big\}\Big)\Big]\\
&=e^{c\operatorname{Vol}_{\RR^d}(A+B(0,R))}
\end{align*}
where $c=c(\lambda, \mu, R)$. 
\subsection{Mat\'ern hardcore processes}
In 1960, Mat\'ern introduced several hardcore models for point processes. These processes are dependent thinnings of PPPs and spread more regularly in space than PPPs. Such models are useful when competition for resources exists (e.g. tree or city locations, see \cite{Matern} and references therein). 

Let $\N$ be a marked PPP of intensity $\lambda$, with independent marks $\{T_x\}_{x\in\N}$ uniformly distributed in $[0,1]$. Then, Mat\'ern I/II hardcore processes (MHP I/II) $\N_{\mbox{I}}$ and $\N_{\mbox{II}}$ are defined, for a given $R>0$ by:
\[\N_{\mbox{I}}:=\big\{x\in\N\,:\,\Vert x-y\Vert>R,  \forall y\in\N\setminus\{x\}\big\},\]
\[\N_{\mbox{II}}:=\big\{x\in\N\,:T_x<T_y,  \forall y\in\N\cap B(x,R)\big\}.\]  
Clearly, MHPs are stationary and have finite range of dependence. Note that $\N_{\mbox{I}}\subset\N_{\mbox{II}}\subset\N$. Hence, the facts that $\N_{\mbox{I}},\,\N_{\mbox{II}}$ are almost surely in general position and that they have almost surely no descending chains are inherited from these properties for PPPs. Moreover, inequalities {\bf (D$_2$)} and {\bf (D$_{3^+}$)} are immediate from those for PPPs, and it suffices to show {\bf (V)} for $\N_{\mbox{I}}$. To this end, first assume there is an integer $n$ such that $L=3nR$ and cut $[0,L]^d$ into $3^dn^d$ disjoint sub-boxes $b_j$ of side $R$. Note that if there is a sub-box $b_j$ with $\# (b_j\cap\N)=1$ having all neighboring sub-boxes empty, then $[0,L]^d\cap \N_{\mbox{I}}\neq\emptyset$. Thus,
\[\PP\big[\#\big([0,L]^d\cap\N_{\mbox{I}}\big)=0\big]\leq \PP\Big[\bigcap_j A_j\Big],\]
where $A_j$ stands for the event: `\emph{$\# (b_j\cap\N)\neq 1$ or a neighbor of $b_j$ contains at least a point of $\N$}'. One can choose a collection of $n^d$ sub-boxes $b_j$ so that, if $i\neq j$, events $A_i$ and $A_j$ are independent. So, there exists a constant $c_{13}>0$ such that:
\[\PP\big[\#\big([0,L]^d\cap\N_{\mbox{I}}\big)=0\big]\leq \PP\big[A_{j_0}\big]^{n^d}=e^{-c_{13}n^d}=e^{-\frac{c_{13}}{3^d}L^d}.\]
For general $L\geq 3R$, it suffices to fix $n$ so that $3nR\leq L\leq 3(n+1)R$ to obtain:
\[\PP\big[\#\big([0,L]^d\cap\N_{\mbox{I}}\big)=0\big]\leq \PP\big[\#\big([0,3nR]^d\cap\N\big)=0\big]\leq e^{-\frac{c_{13}}{6^d}L^d}.\]
\subsection{Determinantal processes}
In 1975, Macchi \cite{Macchi} introduced determinantal point processes (DPPs) in order to model fermions in quantum mechanics. These processes also arise in various other settings such as eigenvalues of random matrices, random spanning trees, carries processes when adding a list of random numbers; see for example \cite{Burton, Diaconis, Ginibre, Shirai1, Shirai2}.

Let $K$ be a self-adjoint non-negative locally trace class operator acting on $L^2(\RR^d)$ with integral kernel $k$. For a bounded Borel set $A\subset\RR^d$, let us denote by $P_A$ the projection operator from $L^2(\RR^d)$ onto $L^2(A)$ and by $K_{|A}:=P_AKP_A$ the restriction of $K$ onto $L^2(A)$. The integral kernel $k$ can be properly chosen such that it satisfies the so called \emph{local trace formula}:  
\begin{equation}
\operatorname{tr}K_{|A}=\int_Ak(x,x)\d x\label{traceformula}
\end{equation}
for all bounded Borel set $A\subset\RR^d$. A DPP with kernel $k$ is a simple point process $\N$ whose correlation functions $\rho_m$ satisfy:
\[\rho_m (x_1,\dots,x_m)=\det \big(k(x_i,x_j)\big)_{1\leq i,j\leq m} \]
for all $m\geq 1$ and all $x_1, \dots, x_m\in \RR^d.$ If $\N$ is stationary,  one can fix $K_0:=k(0,0)=k(x,x)>0$ for almost all $x$. Thus, formula (\ref{traceformula}) reduces in this case to:
\begin{equation}
\operatorname{tr}K_{|A}=K_0\operatorname{Vol}_{\RR^d}(A),\label{traceformula2}
\end{equation}
where $\operatorname{Vol}_{\RR^d}(A)$ denotes the volume of $A$ for the Lebesgue measure on $\RR^d$.
We refer to \cite{BHKPV}, \cite{Soshnikov00} and the appendix in \cite{GeorgiiYoo} for more details.      

Theorem 7 in \cite{BHKPV} provides the following useful fact for DPPs: the number of points of a DPP $\N$ with kernel $k$ falling in a bounded set $A$ has the law of the sum of independent Bernoulli random variables with parameters the eigenvalues of the restriction $K_{|A}$ of $K$ to $A$. Thus writing $\lambda_j$ for the eigenvalues of $K_{|A}$, using that for $s\geq 0,\,1-s\leq e^{-s}$ and formula (\ref{traceformula2}), we have:
\begin{align*}
\PP[\# (A\cap\N)=0 ]&=\prod_j(1-\lambda_j)\\
&\leq \prod_je^{-\lambda_j}=e^{-\operatorname{tr}(K_{|A})}\\
&=e^{-K_0\operatorname{Vol}_{\RR^d}(A)}\\
\end{align*}
which implies {\bf (V)}. Similarly, using that $1+s\leq e^{s}$ and formula (\ref{traceformula2}):
\begin{align*}
\EE\left[\left(\frac{3}{2}\right)^{\#(\N\cap A)}\right]&= \prod_j\left(1+\frac{\lambda_j}{2}\right)\\
&\leq \prod_j e^\frac{\lambda_j}{2}=e^{\frac{1}{2}\operatorname{tr}(K_{|A})}\\
&=e^{\frac{1}{2}K_0\operatorname{Vol}_{\RR^d}(A)}.
\end{align*}
This implies {\bf (D$_2$)} using exponential Markov inequality. Actually, one can obtain {\bf (D$_{3+}$)} in a similar way, but we are unable to check assumption (3) in Criterion \ref{crittrans} because of the lack of independence in this case.  

\subsection*{Acknowledgements}
The author thanks Jean-Baptiste Bardet and Pierre Calka, his PhD advisors, for introducing him to this subject and for helpful discussions, comments and suggestions. The author also thanks anonymous referees for their careful reading and for their comments that significantly improve the paper. This work was partially supported by the French ANR grant PRESAGE (ANR-11-BS02-003) and the French research group GeoSto (CNRS-GDR3477).
\bibliography{biblio}
\bibliographystyle{alpha}
\end{document}